\documentclass[12pt]{amsart}

\usepackage{comment,hyperref,color}

\usepackage{breakurl}   

\usepackage{enumitem} 

\usepackage{comment} 

\usepackage[dvips]{graphicx} 

\usepackage{setspace}  

\numberwithin{equation}{section}

\newcommand{\N}{\mathbb{N}}

\newcommand{\Q}{\mathbb{Q}}
\newcommand{\R}{\mathbb{R}}

\newcommand\astr{{{}^\ast\hspace*{-.5pt}\R}}

\newcommand{\ns}[1]{{}^\ast\hspace{-1pt}#1}

\newcommand\astf{{{}^{\ast}\hspace*{-2.5pt}f}}

\title[Infinite lotteries, spinners, applicability of hyperreals]
{Infinite lotteries, spinners, and the applicability of hyperreals}

\author{Emanuele Bottazzi}
\address{E.~Bottazzi\\
	Department of Civil Engineering and Architecture\\
	University of Pavia\\
	Via Adolfo Ferrata 3, 27100 Pavia, Italy, orcid 0000-0001-9680-9549}
\email{emanuele.bottazzi@unipv.it, emanuele.bottazzi.phd@gmail.com}

\author{Mikhail G. Katz}
\address{M.~Katz\\
	Department of Mathematics\\
	Bar Ilan University\\
	Ramat Gan 5290002 Israel, orcid 0000-0002-3489-0158}
\email{katzmik@macs.biu.ac.il}

\date{\today}

\begin{document}


\begin{abstract}
We analyze recent criticisms of the use of hyperreal probabilities as
expressed by Pruss, Easwaran, Parker, and William\-son.  We show that
the alleged arbitrariness of hyperreal fields can be avoided by
working in the Kanovei--Shelah model or in saturated models.  We argue
that some of the objections to hyperreal probabilities arise from
hidden biases that favor Archi\-medean models.  We discuss the
advantage of the hyperreals over transferless fields with
infinitesimals.  In \cite{BK20b} we analyze two underdetermination
theorems by Pruss and show that they hinge upon parasitic external
hyperreal-valued measures, whereas internal hyperfinite measures are
not underdetermined.
\end{abstract}



\maketitle
\tableofcontents

\section{Introduction}  

Since Abraham Robinson introduced his framework for infinitesimal
analysis in the 1960s (see \cite{Ro61} and \cite{Ro66}), a sizeable
literature has developed in connection with the applicability of said
infinitesimal analysis in probability theory, physics, and other
fields.  Despite the growing body of literature featuring such
applications,%
\footnote{\label{f1}Applications to physics (Albeverio et
al.\;(\cite{Al86}, 1986), Faris (\cite{Fa06}, 2006), Van den Berg and
Neves (\cite{Va07}, 2007), Loeb and Wolff (\cite{Lo15}, 2015)), to
probability theory (Nelson \cite{radically}, 1987), to stochastic
analysis (Capi\'nski and Cutland \cite{Ca95}, 1995), to canards
(Diener and Diener \cite{Di95}, 1995), to mathematical economics and
theoretical ecology (Campillo and Lobry \cite{Ca12}, 2012), to error
analysis (Dinis and van den Berg \cite{Di19}, 2019), to Markov
processes (Duanmu et al.\;\cite{Du21}, 2021).}
the recent years have seen a vigorous debate concerning the
applicability of Robinson's framework in the sciences, with a number
of advocates and also a number of critics.  The latter include
Easwaran, Elga, Parker, Pruss, Towsner, and Williamson.  Recent
additions to the literature are Easwaran--Towsner (\cite{et}, 2018),
Pruss (\cite{Pr18}, 2018), and Parker (\cite{Pa19}, 2019).  Easwaran
and Towsner call into question the applicability of Robinson's
framework to the description of physical phenomena.  A rebuttal
appears in Bottazzi et al.\;(\cite{19c}, 2019).  The present article
focuses mainly on the critiques as formulated by Parker, Pruss, and
Williamson.  These authors have questioned the applicability of
hyperreal models in probability.

In the present text and in the sequel article \cite{BK20b}, we analyze
a claim by Alexander Pruss (AP) that hyperreal models are
\emph{underdetermined}, in the sense that, given a model, allegedly
``there is no rational reason to choose a particular infinitesimal
member of an extension to be a value for the probability''
(\cite[Section\;3.1]{Pr18}) of a single event.  To buttress his claim,
AP exhibits measures assigning a different infinitesimal value to the
event.  We argue, however, that all of AP's additional measures are
parasitic in the sense of Clendinnen \cite{Cl89}.%
\footnote{Clendinnen points out the possibility that ``All members of
any set of empirically equivalent but logically distinct theories
might be parasitic on one key theory.  That is each of the other
members of the set might only be able to be formulated by utilizing
the formulation of the key theory.  If this situation obtained,
[differential underdetermination] would not hold; for the predictions
which could be made by using any one of the set of theories would
nevertheless require the selection of the single key theory. So the
making of these predictions would depend on the selection of a unique
theory'' \cite[p.\;76]{Cl89}.}
More specifically, AP ignores a key property of entities such as
functions and measures in Robinson's framework, namely the property of
being \emph{internal}.  The importance of internality stems from the
fact that Robinson's transfer principle only applies to internal
entities.  Meanwhile, we prove that all of AP's additional measures
are \emph{external}.%
\footnote{For a discussion of the notions of internal and external
entities see Section~\ref{f3}.}
Thus, if one considers only internal hyperfinite measures, no
underdetermination occurs.

AP also claims that certain transferless ordered fields properly
extending the reals, such as the Levi-Civita field or the surreal
numbers, may have advantages over hyperreal fields in probabilistic
modeling.  However, we show in \cite{BK20b} that probabilities
developed over such fields are less expressive than real-valued
probabilities, and inferior to probabilities developed over hyperreal
fields.

In more detail, AP (\cite{Pr18}, 2018) attacks Robinson's framework
for mathematics with infinitesimals, claiming that its applications in
probability cannot have any physical meaning.  AP's critique is more
sophisticated than that of Easwaran--Towsner, in that he acknowledges
at the outset that some commonly voiced objections to hyperreal
numbers are unconvincing.%
\footnote{See Sections~\ref{s21} and \ref{s22b} for technical
details.}
Nevertheless, AP claims that hyperreal probabilities are
\emph{underdetermined}, namely that there is no rational reason to
assign a particular infinitesimal probability to non-empty events that
classically have probability zero.  His argument hinges upon the
following:
\begin{itemize}
\item
examples of uniform processes that allegedly do not allow for a
uniquely defined infinitesimal probability for singletons, and
\item
a pair of theorems asserting that for every hyperreal-valued
probability measure%
\footnote{Throughout the paper probability measures are assumed to be
\emph{finitely additive}.  Notice that finitely additive probability
measures include also~$\sigma$\emph{-additive} probability measures,
but some infinite sample spaces admit finitely additive probability
measures that are not~$\sigma$-additive.}
there exist uncountably many others that induce
the same decision-theoretic preferences.
\end{itemize}
We will argue that the \emph{underdetermination} claim is baseless, by
addressing each of these critiques.  The first critique is addressed
in Section~\ref{sec ex intro} of the present article, whereas his pair
of theorems are analyzed in detail in the sequel article \cite{BK20b}.

We note that the \emph{underdetermination} attack is different from a
previous attack against hyperreal probabilities developed by AP in
(\cite{pruss14}, 2014).  In that paper, AP sought to argue that
infinitesimals are
\begin{quote}
	``too small'' to give plausible probabilities of individual
	outcomes in a countably infinite lottery.%
\footnote{Ironically, Reeder has criticized hyperreal infinitesimals
	for allegedly being \emph{too\;big} in (\cite{reeder}, 2017).
	Reeder's claim is refuted by Bottazzi (\cite{Bo19}, 2019).}
\cite[p.\;1052]{pruss14}
\end{quote}
The 2014 attack was countered by Benci et al.\ (\cite{benci2018},
2018, Section\;4.5, pp.\;531--534).  The more recent Prussian charge is
unrelated to the ``smallness'' argument of~\cite{pruss14}.  Instead,
AP argues that it is impossible uniquely to assign an infinitesimal as
the probability of an event.  A similarity between \cite{pruss14} and
\cite{Pr18} is that in both texts AP fails to take into account the
crucial distinction between \emph{internal} and \emph{external} sets
and functions in Robinson's framework.  Meanwhile, the analysis in
Benci et al.\;\cite{benci2018} is insufficient to address the
\emph{underdetermination} claim.

AP asserts that such underdetermination is a feature of infinitesimal
probabilities generally, but his examples and theorems mainly focus on
hyperreal-valued probability functions.  Moreover, he suggests that
other non-Archimedean extensions of the real field could be more
suitable for the development of infinitesimal probabilities.
Significantly, he makes no attempt to present a model of his uniform
processes in such alternative frameworks with infinitesimals.  We will
address AP's claim in Section~\ref{sec intro prusstheorems} and in
\cite[Section\;4]{BK20b}.

In Section \ref{s22b} we point out some common hidden assumptions in
mathematical modeling of physical phenomena, and analyze some common
biases against Robinson's framework.  Such biases include the
following: the insistence on the use of the natural numbers as the
only possible model for the time scale of processes that ``go on for
the rest of time'' (Section\;\ref{sec infinity}) and the claim that
uncountably many hypernatural numbers are not suitable for the
representation of such discrete processes (Section\;\ref{sec
cardinality}).

In Section \ref{sec intro prusstheorems} we highlight the significance
of the transfer principle of Robinson's framework.  AP suggests that
measures taking values in non-Archimedean fields other than hyperreal
fields may be more suitable for the development of infinitesimal
probabilities.  We show that this suggestion overlooks the
significance of the transfer principle.

A common flaw of criticisms of infinitesimal models, as pursued by
opponents of Robinson's framework, is the assumption that certain
properties of the Archimedean accounts for infinite processes must
also be satisfied by every non-Archimedean probability that represents
it.  However this assumption is unjustifiable; see
\cite[Section\;2]{BK20b}.  Moreover, in \cite{BK20b} we show that
there are appropriate hyperfinite representations of the process that
are not underdetermined.

The issues with the article by AP could be classified along the
following lines:
\begin{enumerate}
[label={(P\theenumi)}]
\item
\label{i1}
(philosophical) (a) AP naively assumes that hyperreal models must
mimick the properties of Archimedean ones; see Sections~\ref{sec
infinity}, \ref{sec cardinality}, and \cite[Section\;2]{BK20b}.  (b)
AP fails to establish the relevance of the parasitic external measures
he introduces to buttress his underdetermination charge; see
\cite{BK20b}, Section\;3.4.
\item
\label{i2}
(historical) AP ignores the Klein--Fraenkel criteria for the
utility of a theory of infinitesimals; see Section \ref{sec KF}.
\item
\label{i3}
(consistency) While at the outset AP admits that the
\emph{arbitrariness} claim against the hyperreals is mathematically
incoherent, he lapses into it later in his article; see
Section~\ref{s34}.
\item
\label{i4}
(mathematical) AP makes inappropriate choices in hyperreal modeling.
Adopting more appropriate choices dissolves AP's argument against
hyperreal modeling; see \cite[Section 2.4]{BK20b}.
\end{enumerate}

In the present article we address items \ref{i1}, \ref{i2}, and
\ref{i3}, whereas in \cite{BK20b} we address items \ref{i1} and
\ref{i4}.

\section
{On mathematical representation of physical processes} 
\label{sec ex intro}

AP opens his analysis by articulating some ``intuitions'' based on
scenarios that can be represented by both Archimedean and
non-Archime\-dean probability measures (before going on to discuss his
\emph{underdetermination} theorems that we will analyze in
\cite[Section\;3]{BK20b}).
%
%
Such scenarios can be grouped into the following categories:
\begin{itemize}
\item some infinite processes, such as coin tosses%
\footnote{The example of an infinite collection of coin tosses has
already been used to attack hyperreal-valued probability measures by
Williamson \cite{williamson}, Easwaran \cite{easwaran14}, Parker
\cite{Pa19} and other authors.  Rebuttals of this argument can be
found for instance in Weintraub (\cite{We08}, 2008), Bascelli et
al.\;(\cite{14a}, 2014), Hofweber (\cite{Ho14} 2014), and Howson
(\cite{howson19}, 2019).}
\cite[Sections\;3.2 and 4.2]{Pr18}, or the estimate of the value of a
utility one can have ``every day for eternity''
\cite[Section\;3.5]{Pr18};
\item a pair of uniform processes over a single sample space,
exemplified by the motion of two spinners (rotating pointers)
\cite[Sections\;3.2 and 3.3]{Pr18} and by a pair of uniform lotteries
over~$\N$ \cite[Section\;4.1]{Pr18}.
\end{itemize}

We will discuss the details of AP's representation of the two spinners
in \cite[Section\;2]{BK20b}.  Here we will comment more generally on
the issue of mathematical representation of physical processes.  We
argue that certain physical processes may admit distinct mathematical
models.  For such processes, there does not exist a unique,
well-defined model that would represent them in a way resembling
anything like an isomorphism.  Such a perspective is accepted even by
mathematicians and philosophers who adopt a responsible variety of
mathematical realism, as discussed in Section~\ref{sec infinity} (see
also Bottazzi et al.\;\cite{19c}, 2019, Section\;1).  We further argue
that hyperreal models can be used on par with Archimedean models based
upon the Cantor--Dedekind representations of the continuum.  This
issue is also dealt with by Herzberg (\cite{He07}, 2007, Section\;4).

As acknowledged by AP, a common objection to the use of hyperreal
fields in such representations, namely that such fields are
arbitrarily specified, can be countered in several ways.  Such an
objection can be countered either by working in the Kanovei--Shelah
definable hyperreal field \cite{KS} or by using suitable saturated
hyperreal fields; see Section~\ref{s21}.

In Sections~\ref{sec infinity} and \ref{sec cardinality} we address
further hidden assumptions in Archimedean mathematical descriptions of
the infinite processes represented by coin tosses from the
aforementioned point of view that rejects the postulation of a unique
mathematical representation.  In particular, we argue that some
commonly voiced objections to the use of hyperreals in the
representation of scenarios involving an infinity of events stem from
such hidden assumptions concerning mathematical representation of
physical events.

\subsection{Are hyperreal fields arbitrary?}
\label{s21}

AP acknowledges at the outset the failure of the commonly voiced
objection of \emph{arbitrariness} (and even that of
\emph{ineffability}), namely that one cannot specify a particular
hyperreal extension of~$\R$ (see Section~\ref{s22b} for details).
Such an objection was voiced by Alain Connes%
\footnote{See (\cite{connes}, 2004, p.\;14) where Connes describes
Robinson's framework as``some sort of chimera.''}
and others.  The objection is specifically refuted by the
Kanovei--Shelah definable hyperreal field \cite{KS}; see also
\cite{herzberg}.  For a rebuttal of the Bishop--Connes critique see
Katz--Leichtnam (\cite{13d}, 2013), Kanovei et al.\;(\cite{15a},
2015), and Sanders (\cite{Sa20}, 2020).

Furthermore, AP acknowledges that the \emph{arbitrariness} objection
can also be refuted by working with a hyperreal field defined up to an
isomorphism, and that for suitable cardinals~$\kappa$, there is a
unique-up-to-isomor\-phism $\kappa$-saturated hyperreal field of
cardinality~$\kappa$.%
\footnote{\label{f7}More precisely, the condition is that an infinite
cardinal~$\kappa$ should either be inaccessible or satisfy~$2^\kappa=
\kappa^+$ (that is, the continuum hypothesis holds at~$\kappa$).  For
details see Keisler (\cite{Ke94}, 1994, Section\;11) and further
references therein.}
We will source such ``concessions'' by AP in Section~\ref{s22b}.

\subsection
{Pruss admits failure of arbitrariness/ineffability claims}
\label{s22b}

AP mentions a worry that
\begin{quote}
the choice of a hyperreal extension appears to be not only arbitrary
but ineffable {\ldots} -- we cannot successfully refer to a particular
extension, and so a particular extension cannot reflect our credences
{\ldots} \cite[Section\;1]{Pr18}
\end{quote}
However, he immediately acknowledges that ``the ineffability argument
does not apply to all extensions of the reals, and even as restricted
to the hyperreals it is \emph{unsuccessful}'' (ibid.; emphasis added).

Thus AP acknowledges at the outset that the commonly voiced
objection of \emph{arbitrariness} (and even that of
\emph{ineffability}), namely the claim that one cannot specify a
particular hyperreal extension of~$\R$, is unsuccessful, and
specifically refuted by the Kanovei--Shelah definable hyperreal field
(see \cite{KS}, \cite{herzberg}):
\begin{quote}
[B]y leveraging the idea that even when it is difficult to specify a
particular ultrafilter, one can specify sets of ultrafilters, Kanovei
and Shelah (2004) explicitly defined a particular free ultrafilter on
a particular infinite set%
\footnote{This characterisation by AP of the Kanovei--Shelah
technique contains a mathematical inaccuracy.  The technique does not
exploit an ultrafilter on an infinite set.  Rather, it exploits a
maximal filter \emph{in a particular algebra of subsets} of a certain
infinite set (the algebra in question does not contain all subsets!).}
%
%
	%
%
{\ldots} Furthermore, Kanovei and Shelah used their construction to
make an explicitly specified extension of the reals (an iteration of
the hyperreal extension using this ultrafilter) having further
desirable properties.  \cite[Section~2]{Pr18}
\end{quote}
Moreover, AP acknowledges that the \emph{arbitrariness} objection
can also be refuted by working with a hyperreal field defined up to an
isomorphism:
\begin{quote}
[W]e can specify a set of hyperreals up to isomorphism.  For some
cardinals~$\kappa$, there is a unique-up-to-isomor\-phism
$\kappa$-saturated non-standard real line of cardinality~$\kappa$
{\dots} And there might be some non-arbitrary way to choose the
cardinal~$\kappa$, perhaps a way matching the particular problem under
discussion.  \cite[Section 2]{Pr18}
\end{quote}

\subsection{Shift-invariance hypothesis} 
\label{sec infinity}

Various scenarios involving infinite processes have been discussed by
AP and other authors, including Easwa\-ran, Parker, and William\-son.
Such discussions often exhibit a bias in favor of Archimedean models,
which feature
\begin{itemize}
	\item a countable infinity of events, and
	\item events that are ordered in time (rather than simultaneous).
\end{itemize}

A typical process that is modeled with such hidden assumptions is the
outcome of an infinite amount of coin tosses.  A number of arguments
against hyperreal probabilities for infinite coin tosses hinge upon
the events~$H(n)$ that AP defines as follows:
\begin{quote}
starting with day~$n$, it's all heads for the rest of time.
\cite[Section\;3.2]{Pr18}
\end{quote}
AP models such a process by a sequence of tosses labeled by $\N$, and
makes the following assumption:
\begin{quote}
\textbf{Shift-invariance hypothesis}: Events~$H(n)$ and~$H(m)$ are
isomorphic for all~$m,n \in \N$
\end{quote}
(the term \emph{shift-invariance hypothesis} is ours).  Meanwhile,
alternative models of these infinite processes can be obtained with
hyperfinite techniques, as discussed for instance by Benci et al.\
\cite[pp.\;44--46]{bbd2}; see also Nelson (\cite{radically}, 1987) and
Albeverio et al.\;(\cite{Al86}, 1986).  These models show that the
shift-invariance hypothesis is spurious, since it does not hold in a
hyperfinite representation of the infinite collection of coin tosses.%
\footnote{Howson argued that the events~$H(n)$ and~$H(m)$ are not
equiprobable whenever~$n\ne m$ even in the Archimedean model where the
sample space is the~$\sigma$-algebra generated by the cylinder sets
in~$\{0,1\}^{\N}$ (\cite{howson19}, 2019, Section\;3).  A similar
observation was made by Benci et al.\;\cite[pp.\;21--22]{benci2018}.
Both arguments have been addressed by Parker \cite{Pa19}.}

The shift-invariance hypothesis is often justified by an appeal to the
``physical structure'' of the infinite process.  Thus, Williamson
writes: ``A fair coin will be tossed infinitely many times at one
second intervals'' in (\cite{williamson}, 2007) on page 174.  By the
middle of page 175, he is ready to claim an ability to
\begin{quote}
map the constituent single-toss events of~$H(1\ldots)$ one-one onto
the constituent single-toss events of~$H(2\ldots)$ in a natural way
that preserves the \emph{physical structure} of the set-up just by
mapping each toss to its successor. (ibid., p.\;175; emphasis added)
\end{quote}
Williamson appears to be taking for granted a ``physical structure''
possessing a considerable supply of physical seconds.

The same assumption appears in the more recent text by Parker
(\cite{Pa19}, 2019).  Parker implicitly assumes that a countable
sequence of coin tosses is physically feasible, and bases his
\emph{Isomorphism Principle} \cite[p.\;4]{Pa19} on such an assumption.

Bascelli et al.\;(\cite{14a}, 2014) analyzed similar biases in favor
of modeling based upon a countable infinity of time-ordered events in
Easwaran (\cite{easwaran14}, 2014).  Namely, an assumption of a
countable time-ordered model already involves a full-fledged
idealisation lacking a referent.

What Williamson and Parker fail to recognize is that, even from the
viewpoint of a responsible variety of mathematical realism, a
mathematical description of a physical event typically involves some
level of idealisation and introduces some spurious properties (as
already argued for instance in \cite[Section 1]{19c}).  As a
consequence of such idealisation, it is not possible to claim that a
physical process has a unique well-defined sample space, or that one
particular sample space provides the only correct mathematical
description of a physical process. For more details, we refer to
\cite[Section\;1.4]{19c} and to references therein.

In the coinflip case, it is obvious that it is physically impossible
to flip a coin, as Parker would have it, ``infinitely many times, at
times~$t_0 + n$ seconds for~$n = 0, 1, 2,\ldots$'' \cite[p.\;8]{Pa19}.
Nevertheless, it is possible to \emph{model} this situation as a
sequence of coin tosses over~$\N$, or with other notions of number, as
already mentioned.

Furthermore, starting with a physical intuition of something going on
``for the rest of time,'' there are several possibilities of
formalizing such an intuition mathematically.  One way is to interpret
time increments as ranging over the traditional~$\N$.  An alternative
way of modeling such an intuition would be to postulate that the
``time'' in question comes to an end rather than goes on indefinitely,
given the likelihood of physical armageddon expected by some modern
theories in astrophysics.  If so, then finite and hyperfinite
modeling, which postulate such a final moment, are arguably more
faithful to physical intuition than modeling by~$\N$.  In this sense,
an assumption that an intuition of ``for the rest of time''
necessarily refers to~$\N$, involves circular reasoning, as the
conclusion is built into the premise.  Attempting to base ``intuitive
reasons" against infinitesimal probabilities on such an idealized
model, as AP does in \cite[Section\;3.2]{Pr18}, begs the question as
to why one assumes precisely such a model rather than, say, a
hyperfinite number of simultaneous coin flips.%
\footnote{Or finite nonstandard number of flips in Nelson's framework;
see Section~\ref{sec def transfer}.}
Indeed, the model chosen by Easwaran, Parker, Pruss, and Williamson
predetermines the outcome of their analyses.  The shift-invariance
hypothesis is analyzed further in Section~\ref{s25}.

\subsection{Cardinality objection}   
\label{sec cardinality}

AP puts forth the following objection to the use of hypernatural
numbers:
\begin{quote}
	[I]t turns out that for any positive infinite number~$M$
	in~$\astr$, if~$\astr$ has a collection of hypernaturals, then
	there will be uncountably many (in external cardinality)%
\footnote{AP's parenthetical comment referring to \emph{external
cardinality} indicates that he is aware of the distinction between
internal and external entities (in this case, cardinality).  Six years
prior to the online publication of \cite{Pr18}, he referred to
\emph{internal cardinality} in his posting \cite{Pr12}.  However, AP
tails to take into account the distinction between internal and
external hyperreal probabilities, as we will show in \cite{BK20b},
Section 2.5.}
hypernaturals between~$1$ and~$M$ (Pruss 2014, Appendix).  And so the
countable number of future days that we've imagined is not what is
counted by~$M$: instead,~$M$ counts the number of members of the
uncountable set
$\{1,2,\ldots,M\}$ of hypernaturals.  (Pruss \cite{Pr18}, 2018,
Section\;3.5)
\end{quote}
What AP is claiming is that if one is interested in countably many
``future days'' (i.e., trials), hypernaturals do not provide an
accurate model by cardinality considerations.  However, his
cardinality objection is not valid, for the following reason.  Skolem
\cite{skolem} already developed elementary (in the sense of PA)
extensions~$\N_{Sk}$ of~$\N$ in the 1930s.  Skolem's precedent was
clearly acknowledged by Robinson (\cite{Ro66}, 1966, pp.\;vii, 88,
278), who noted that ``Skolem's method foreshadows the ultrapower
construction'' (op.\;cit., p.\;88).  Being built out of equivalence
classes of (definable) sequences of natural numbers,~$\N_{Sk}$
naturally embeds in~$\astr$ (for details see Kanovei et
al.\;\cite{kanovei x skolem}, 2013).  If one's interest is in
countable structures only, one can proceed as follows:
\begin{enumerate}
	\item
	construct a countable extension~$\N_{Sk}$ of~$\N$ following
	Skolem;
	\item
	form the field of fractions~$F$ of~$\N_{Sk}$;
	\item~$F$ is then an ordered field properly including~$\Q$.
\end{enumerate}
In particular, there will be only countably many numbers in such an
extension~$F$, and hence countably many numbers in the set
$\{1,2,\ldots,M\}$.  An identical rebuttal applies to AP's rejection
of a nonstandard solution to the paradox of Thomson's lamp in
(\cite{Pr18b}, 2018, p.\;41).

\subsection
{Standard model of the naturals and coinflipping}
\label{s25}

In this section we will examine the relation of the so-called standard
model of arithmetic to modeling infinite processes such as infinite
lotteries, coin flips, etc.  It is possible to disassociate the issue
of scientific modeling (in physics, probability, etc.) from the issue
of putative existence of a standard model (a.k.a.\;the intended
interpretation) of~$\N$.  Even modulo such an~$\N$ along the
Cantor--Dedekind lines (not along the Nelson lines), one can question
the Pruss--Williamson (PW) assumption that~$\N$ can be embedded in
physical time.  PW make no effort to justify the assumption, which is
surprising for publications in venues such as \emph{Analysis} and
\emph{Synthese}.

In the following, we adopt the analysis of Kuhlemann (\cite{Ku18},
2018).  When PW speak of performing a trial every second (or day) from
now to eternity and of the ``physical structure'' of the process, they
may be referring to either metalanguage natural numbers or the object
language natural numbers.

Logicians make a distinction between, on the one hand, metalanguage
naturals, and, on the other, the object language naturals (e.g.,
numbers in the putative standard model a.k.a.\;intended
interpretation).  Thus, Simpson denotes metanaturals by $\omega$ to
distinguish them from $\N$ \cite[pp.\;9--10]{Si09}.  At best,
metanaturals can be related to as a sorites-type subcollection (of the
object language naturals $\N$) which does not exist as a set, blocking
implementation of a set of trials indexed by metanaturals (see
\cite[p.\;255]{Di19} for a related model of the sorites paradox).

If PW mean to refer to metanaturals, the analysis above would
undermine the shift-invariance hypothesis (see Section~\ref{sec
infinity}) and the PW attempt to identify $H(1)$ with~$H(2)$
``naturally''.%
\footnote{Williamson actually uses the term \emph{natural} in
reference to applying the shift to physical processes.  AP actually
speaks of ``[t]he countable number of future days''
\cite[Section\;3.5]{Pr18}.}

If, on the other hand, PW are referring to the object language
naturals, then they are already making an assumption favoring one type
of idealisation over another, so that their conclusion is built into
their premise.  $\N$ is not naturally built into intuitions and
thought experiments involving lotteries and coinflips (though it may
be built into the type of undergraduate mathematical training that PW
received).

\section{On the strength of theories with infinitesimals} 
\label{sec intro prusstheorems}


We will define the notions of internal and external objects in
Section\;\ref{f3}, and present the transfer principle of Robinson's
framework in Section\;\ref{sec def transfer}.  These notions play a
major role in mathematical modeling with hyperreal numbers.  The
significance of the transfer principle is also related to the
historical development of mathematical theories with infinitesimals,
as discussed in Section\;\ref{sec KF}.  Thus these notions will be
central to our discussion of the infinitesimal models of the uniform
processes proposed by AP (see \cite{BK20b}, Section\;2), and of AP's
pair of theorems (see \cite{BK20b}, Section\;3).

We will also evaluate AP's claim that transferless non-Archimedean
extensions of the real numbers might be more suitable for the
development of infinitesimal probabilities.  In Sections\;\ref{sec KF}
and \ref{s34} we will elaborate on some consequences of the absence of
transfer in the surreal numbers and the Levi-Civita field.  What this
entails for their applicability or otherwise to probability theory is
discussed in detail in \cite{BK20b}, Section 4.

\subsection{Constructing hyperreal fields} 
\label{f3}

It is well known that fields~$\astr$ of hyperreal numbers can be
obtained by the so-called \emph{ultrapower construction}.  In this
approach, one sets~$\astr = \R^\N/\mathcal{U}$, where~$\mathcal{U}$ is a
nonprincipal ultrafilter over~$\N$. The operations and relations on
$\astr$ are defined from the quotient structure.  For instance, given~$x
= [x_n]$ and~$y=[y_n]$, we set~$x+y=[x_n+y_n]$ and~$x \cdot
y=[x_n\cdot y_n]$.  We have~$x < y$ if and only if~$\{n\in\N\colon
x_n<y_n\}\in \mathcal{U}$.

Let~$\mathcal P=\mathcal P(\R)$ be the power set of~$\R$.  Then the
star transform produces the object~${}^{\ast}\mathcal P$.  An
\emph{internal} subset~$A\subseteq\astr$ of~$\astr$ is by definition a
member of~${}^{\ast}\mathcal P$.  More concretely, in the ultrapower
construction an internal subset~$A\subseteq\astr$ is represented by a
sequence~$(A_n)$ of subsets~$A_n\subseteq\R$.  Here an element
$[x_n]\in\astr$ belongs to~$A=[A_n]$ if and only if~$\{n\in\N\colon
x_n\in A_n\} \in \mathcal{U}$.  A subset of~$\astr$ which is not
internal is called \emph{external}.

More generally, in the context of the star transform from the
superstructure over~$\R$ to the superstructure over~$\astr$, a set~$A$
of the latter is internal if and only if it is a member of
${}^\ast\hskip-2pt Z$ for some~$Z$ in the superstructure over~$\R$.
For further properties of the ultrapower construction of hyperreal
numbers and the superstructures, see Fletcher et al.\ (\cite{17f},
2017) and Goldblatt (\cite{goldblatt}, 1998).

\subsection
{The transfer principle of Robinson's framework} 
\label{sec def transfer}

Kanovei et al.\ describe the transfer principle of Robinson's
framework as
\begin{quote}
	a type of theorem that, depending on the context, asserts that
	rules, laws or procedures valid for a certain number system,
	still apply (i.e., are ``transferred") to an extended number
	system. (\cite{18i}, 2018, p.\;113)
\end{quote}
%

The transfer principle asserts that the \emph{internal objects} of
Robinson's framework satisfy all the first-order properties of the
corresponding classical objects.

The simplest examples of transfer involve the extension of sets and
functions via the~${}^\ast$ map.  For instance, a continuous function
$f\colon \R \rightarrow \R$ is a function that satisfies the formula
\begin{align*}
	&\forall x_0 \in \R\ \forall \varepsilon \in\R, \varepsilon > 0\ \exists \delta \in\R, \delta > 0\\
	&\forall x \in \R\ \left( |x-x_0|<\delta \rightarrow |f(x)-f(x_0)|<\varepsilon \right).
\end{align*}
The function~$f$ is extended to a function~$\astf\colon\astr\rightarrow
\astr$ that satisfies
\begin{align*}
&\forall x_0 \in \astr\ \forall \varepsilon \in\astr, \varepsilon > 0\
\exists \delta \in\astr, \delta > 0\\ &\forall x \in \astr\ \left(
|x-x_0|<\delta \rightarrow |\astf(x)-\astf(x_0)|<\varepsilon\right)
%
%
\end{align*}
(an internal function satisfying this formula is sometimes called
$\ns$conti\-nuous).  Moreover,~$\astf$ satisfies the Intermediate
Value Theorem, the Mean Value Theorem, and every other first-order
property of~$f$.

We now turn to properties of sets under extension.  For instance, the
Archimedean property of~$\R$ is expressed by the formula
\begin{equation}
\forall x,y \in \R\ \big((0<x \land x<y) \rightarrow \exists n \in \N
\,(y<nx)\big).
\end{equation}
An application of the transfer principle to the above formula yields
\begin{equation}
\label{e22}
\forall x,y \in \astr\ \big((0<x \land x<y) \rightarrow \exists n \in
\ns{\N} \,(y<nx)\big).
\end{equation}
The latter formula needs to be distinguished from the former, since,
as is well known, a ring extension of~$\R$ with infinitesimal elements
is non-Archimedean.%
\footnote{The fundamental difference between the two formulas is that,
in \eqref{e22}, the variable~$n$ can take infinite hypernatural
values.}

The first-order properties are preserved also by certain sets and
functions that are not of the form~${}^\ast\hspace{-2pt}{X}$ for some
classical~$X$.  A relevant example is given by a \emph{hyperfinite
set}, i.e., a set that can be put in an (internal) one-to-one
correspondence with a set of hypernatural numbers of the form~$\{x \in
\ns{\N}\colon x \leq H\}$.  Hyperfinite sets being internal, the
transfer principle ensures that hyperfinite sets have the same
first-order properties as finite sets.  As a consequence, hyperfinite
sets can be routinely applied to a wide variety of problems.  For a
discussion of selected applications, we refer to Arkeryd et al.\
(\cite{theoryapps}, 1997).
	
Katz--Sherry (\cite{13f}, 2013) suggest that the transfer principle
can be interpreted as a formalisation of the \emph{law of continuity}
of Leibniz; see also Sherry--Katz (\cite{14c}, 2014).  Such a
connection was first mentioned by Robinson (\cite{Ro66}, 1966,
p.\;266).  Robinson's historical chapter 10 has occasioned a
reappraisal of the legacy in infinitesimal analysis of pioneers like
Fermat \cite{18d}, Gregory \cite{18f}, Leibniz \cite{18a}, Euler
\cite{17b}, and Cauchy\;\cite{20a},\;\cite{19a}.

We emphasize that the transfer principle applies only to
\emph{internal} entities of Robinson's framework.  Note that external
entities do not exist in Nelson's framework \emph{Internal Set Theory}
\cite{nelson}.%
\footnote{\label{f14}In Nelson's framework, set theory is enriched by
a one-place predicate \textbf{st}.  The formula~$\textbf{st}(x)$
asserts that an entity~$x$ is standard.  The ZFC axioms are enriched
by the addition of further axioms governing the interaction of the new
predicate with the axioms of traditional set theory.  Infinitesimals,
say~$\epsilon$, are found within the ordinary real line, and satisfy
$0<|\epsilon|<r$ for all standard~$x\in\R^+$.  It is shown in
\cite{nelson} that the new theory is conservative with respect to ZFC.
A related system was developed independently by Hrbacek \cite{Hr78}.
For further details, see Fletcher et al.\;(\cite{17f}, 2017) and
Hrbacek--Katz (\cite{Hr20}, 2020).}
For this reason, attempted arguments from first principles that do not
take into account Nelson's framework are not actually based on first
principles as they are claimed to be, but rather involve an unspoken
commitment to a specific set-theoretic framework (for instance,
Zermelo--Fraenkel set theory plus the Axiom of Choice) expressed in
the~$\in$-language.  This is done at the expense of other possible
foundational frameworks.%
\footnote{More specifically, the unspoken commitment typically
involved is to a set theory in the~$\in$-language rather than a set
theory in the~$\in$-$\mathbf{st}$-language; see note~\ref{f14}.}
Thus, from the point of view of Internal Set Theory, internal
probability measures are no less \emph{underdetermined} than the
traditional ones.

%

\subsection
{Usefulness of infinitesimals: Klein and Fraenkel} 
\label{sec KF}

During the opening decades of the 20th century, both Felix Klein%
\footnote{In recent decades, there has been a deplorable attempt by
Mehrtens (\cite{Me90}, 1990), Gray (\cite{Gr08}, 2008), and others to
discredit Klein both mathematically and politically.  A rebuttal
appears in Bair et al.\;(\cite{18b}, 2017).}
(\cite{Kl08}, 1908) and Abraham Fraenkel (\cite{Fr28}, 1928) formulated
a pair of criteria to gauge the success of theories with
infinitesimals.  These criteria are
\begin{enumerate}
\item
the availability and provability (by infinitesimal techniques) of the
Mean Value Theorem, and
\item
the introduction of the definite integral in terms of infinitesimal
increments.
\end{enumerate}
For a detailed discussion, see Kanovei et al.\ (\cite{18i}, 2018).
Klein and Fraenkel both observed that the infinitesimal theories
available at the time (including the Levi-Civita field) did not enable
a satisfactory treatment of these topics.
%
%
When Robinson introduced his framework for analysis with
infinitesimals, Fraenkel related to Robinson's framework as an
important accomplishment that finally solved the old problem of
developing a usable non-Archimedean field.
%
%
Thanks to the transfer principle, in Robinson's framework it is
possible to prove the Mean Value Theorem (MVT) and to define the
Riemann integral of a continuous function by means of hyperfinite
summation of infinitesimal\;terms.

Notice that the MVT and the definition of an integral, and in general
the development of a calculus on non-Archimedean structures that
extends the real calculus, require an extension of real functions.
Costin et al.\ observe that
\begin{quote}
A longstanding aim has been to develop analysis on [the surreal
numbers] as a powerful extension of ordinary analysis on the
reals. This entails finding a natural way of extending important
functions~$f$ from the reals to the reals to functions~$f^*$ from the
surreals to the surreals, and naturally defining integration on the
$f^*$.  \cite[Abstract]{integral}
\end{quote}
In Robinson's framework, such an extension is provided by the $\ast$
map.  Meanwhile, it is still an open problem to define well-behaved
extensions of real functions to the surreals or to the Levi-Civita
field.  For the surreal numbers, the problem is caused by the
necessity to define functions from the \emph{simplicity hierarchical
structure} of the surreal number tree.%
\footnote{For instance, no surreal extension of the sine function
usable in ordinary mathematics is as yet available; see e.g.,
Kanovei's remark at \url{https://mathoverflow.net/a/307114}}

Meanwhile, every real continuous function admits a canonical extension
to a continuous function on the Levi-Civita field. However, this
extension does not preserve many properties of the original real
continuous functions, such as an Intermediate Value Theorem or a Mean
Value Theorem.%
\footnote{In this context, Bottazzi suggested an analogy between these
extensions and external functions of Robinson's framework
\cite{Bo18}.}

Moreover, the surreal field and the Levi-Civita field satisfy only
some restricted versions of the Klein--Fraenkel criteria.  Namely, in
the surreal numbers it is not possible to prove the MVT, and it is
still an open problem to define an integral (see for instance Costin
et al.\;\cite{integral} and Fornasiero \cite{fornasiero}).  Meanwhile,
for the Levi-Civita field, Shamseddine showed that the MVT is valid
only for analytic functions \cite{mvtcivita}. Consequently, this
theorem fails for instance for the extension of non-analytic real
continuous functions.  The Levi-Civita field does have a notion of
integral in dimensions~$1$,~$2$ and~$3$, but the integral is not
defined in terms of sums of infinitesimal contributions; see
Berz--Shamseddine (\cite{berz+shamseddine2003}, 2003), Shamseddine
(\cite{shamseddine2012}, 2012), Shamseddine--Flynn
(\cite{shamseddine+flynn}, 2018).  In addition, it turns out that the
extensions of continuous but non-analytic real functions are not
measurable%
\footnote{Recall that every real continuous function is
Lebesgue-measurable and, if it is defined over a closed interval, it
also has a well-defined Riemann integral. The failure of measurability
for the extensions of non-analytic real functions is a significant
limitation for the measure theory on the Levi-Civita field and it is
also a blatant failure of transfer for this field.}
(Bottazzi \cite{bottazzi forthcoming}).

When Klein and Fraenkel formulated their criteria for the evaluation
of non-Archimedean extensions of the reals (see \cite{18i}), a number
of such non-Archimedean options were available, including the
Levi-Civita field.  It is those ordered fields that Klein and Fraenkel
were referring to when they expressed disappointment with the (then)
current rate of progress, when even the Mean Value Theorem was not yet
provable using infinitesimal analysis.
%
%
AP's critique of Robinson's framework fails to take into account the
fact that at present, Robinson's is the only theory of infinitesimals
that meets the Klein--Fraenkel criteria of utility.  These criteria
are prerequisites for a measure or probability theory.

\subsection
{Infinitesimals without transfer}    
\label{s34}

The importance of the transfer principle in non-Archimedean extensions
of the real numbers can be better appreciated if one considers what
happens when this principle is not available.  We will refer to such
non-Archimedean fields as \emph{transferless}.

Thus, in Henle's \emph{non-nonstandard analysis} \cite{henle} what is
available is a weak form of transfer that applies only to equations,
inequalities and their conjunctions, but not to their disjunction.  As
a consequence, some properties of ordered fields fail, and Henle's
extension is only a partially ordered ring with zero divisors.%
\footnote{In some cases, working with number systems lacking the
	habitual properties can lead the author into error.  Thus, Laugwitz
	pointed out that Henle's article lapses into using denominators when
	working with a ring.  Laugwitz goes on to invite the reader to
	``rewrite the relevant passages''\;\cite{La99}.}

Similarly, in the Levi-Civita field the absence of a transfer
principle makes it necessary to prove individually many theorems of
the calculus, such as the Intermediate Value Theorem and the Mean
Value Theorem for analytic functions.  For a more detailed discussion,
see Shamseddine--Berz (\cite{shamseddine berz analysis}, 2010) and
Shamseddine (\cite{mvtcivita}, 2011).  In addition, currently it is
possible to extend only analytic real functions to the Levi-Civita
field in a way that preserves their first-order properties, as shown
by Bottazzi (\cite{Bo18}, 2018).

In other transferless fields the situation might be even more
difficult; we are not aware of any research towards establishing some
(even limited) forms of transfer in such settings.  Nevertheless, AP
suggests that there are ``multiple methods'' of developing
infinitesimal probabilities in transferless fields:
\begin{quote}
[H]yperreals are not the only way to get infinitesimals. There are
multiple methods that do not make use of anything like the arbitrary
choice of an ultrafilter.%
\footnote{Here AP seems to lapse into the arbitrariness charge against
Robinson's framework discussed in Section\;\ref{s21}.  Notice that
this passage comes after AP's admission that it is possible to
uniquely specify some particular hyperreal fields.}  %
\cite[Section\;2]{Pr18}
\end{quote}
AP argues that, as a consequence,
\begin{quote}
The friend of infinitesimal probabilities has a real hope of
non-arbitrarily specifying a particular field of infinitesimals.
(ibid.)
\end{quote}
Moreover, he claims that certain transferless fields, namely the
surreal numbers, fields of Laurent series or the Levi-Civita field,
have some advantages over hyperreal fields.  Thus AP writes:
\begin{quote}
[T]he surreals have the advantage of being exhaustively large, large
enough that they escape the cardinality arguments against regularity
of Pruss (2013a).  The fields of formal Laurent series and the
Levi-Civita field, on the other hand, have the advantage of being
elegantly small. (On the other hand, the Kanovei--Shelah field, while
mathematically fascinating, probably has little going for it in this
context.)%
\footnote{AP's parenthetical claim that the KS model has little
advantage over the transferless fields mentioned is questionable,
since the KS model can be used just as well as any
traditional~$\R^{\N}/\mathcal{U}$ model (see Section~\ref{f3}),
particularly with the saturation improvement provided by
Kanovei--Shelah \cite{KS}. In particular, the KS model has the
advantage of the transfer principle.}
\cite[Section\;2]{Pr18}
\end{quote}
It can be argued that the main advantage of the field of Laurent
series and the Levi-Civita field is that of being definable from the
natural numbers in a choice-free manner.  If one works in the von
Neumann--Bernays--G\"odel set theory with global choice, then the
surreal numbers are also uniquely specified.  It is also true that
many hyperreal fields obtained as ultrapowers of~$\R$ are not
definable in a choice-free manner; however, as AP acknowledges,
suitable~$\kappa$-saturated fields of hyperreal numbers are uniquely
specified up to an isomorphism.%
\footnote{See Section~\ref{s21}.}

However, the issue at hand is not whether a non-Archimedean field is
definable from the natural numbers without additional parameters.  The
real issue is the applicability of such a field.  In this regard, the
Levi-Civita field has some limited applications in automatic
derivation \cite{shamseddine berz analysis} and in the description of
physical phenomena \cite{flynn}, while there is a vast literature of
applications of hyperreal fields.%
\footnote{Some relevant examples can be found in note~\ref{f1}.}

That such applications are possible is due in particular to the
transfer principle of Robinson's framework.  Significantly, the
principle is not discussed by AP in \cite{Pr18}.%
\footnote{It is a matter of public record that AP is aware of the
transfer principle, since he mentions it both in \cite[Appendix]{Pr18}
and in \cite[p.\;41]{Pr18b}.}
Instead, AP expresses enthusiasm about the surreal numbers, Laurent
series, and Levi-Civita fields, but fails to explain their
shortcoming, namely lack of transfer.  Significantly, AP does not
develop a model either for the infinite coin tosses or for his
spinners in any of these transferless structures.  Thus, his claim
that these non-Archimedean fields may be suitable for the development
of an infinitesimal probability is baseless.

Indeed, in \cite[Section\;4]{BK20b} we argue that attempts to develop
infinitesimal probabilities over the surreal numbers or over the
Levi-Civita field encounter a number of difficulties.

For instance, in order to accomplish anything with the surreals one
would have first to import the transfer principle via an
identification of maximal class-size surreals with maximal class-size
hyperreals, as mentioned by Ehrlich \cite[Theorem 20]{ehrlich}.
Without such an identification and without the transfer principle of
Robinson's framework, it is currently not possible to develop a
measure theory on the surreal numbers (for more details, see
\cite[Section 4.1]{BK20b}).%
\footnote{It should be noted that the omnific surnaturals do not
satisfy the axioms of Peano Arithmetic; e.g., there exist surnaturals
$p,q$ such that $p^2=2q^2$.  For further details see \cite{Mo93},
\cite{Ru14}, \cite{Je13}.}

The Levi-Civita field does obey a type of transfer principle, albeit
limited to a particular extension of real analytic functions; see
Bottazzi (\cite{Bo18}, 2018).

One could also work directly with the measure-theoretic tools
available in the various theories to define non-Archimedean
probability measures; however, there seem to be some difficulties.

Consider for instance the case of the Levi-Civita field, where a
uniform measure is currently under development by Shamseddine and Berz
\cite{berz+shamseddine2003}, Shamseddine \cite{shamseddine2012},
Shamseddine and Flynn \cite{shamseddine+flynn} and Bottazzi
\cite{bottazzi forthcoming}.  So far this uniform measure is not able
to accommodate more than locally analytic probability functions, and
has no notion of hyperfiniteness comparable to that of Robinson's
framework.  This example is discussed further in \cite{BK20b},
Section\;4.2.

\section{Conclusion}

We have examined some commonly voiced objections to the use of
hyperreal numbers in mathematical modeling.  Many of these objections
are based upon naive assumptions regarding the possibility of uniquely
specifying some hyperreal fields, and upon examples of infinite
processes that allegedly do not allow for a uniquely defined
infinitesimal probability for singletons.

The first objection, namely that it is allegedly not possible to
specify a hyperreal field in a unique way, is refuted by the
Kanovei--Shelah definable hyperreal field and by the fact that, for
suitable infinite cardinals~$\kappa$, there is a
unique-up-to-isomorphism~$\kappa$-saturated hyperreal field of
cardinality~$\kappa$.  Note that this rebuttal is accepted also by
some detractors of Robinson's framework for analysis with
infinitesimals.

With regard to objections based upon the analysis of certain infinite
processes, we have observed that physical processes often admit
distinct mathematical models, so that there does not exist a unique,
well-defined model that would represent them in a way resembling
anything like an isomorphism.  Thus we have shown that, for some
commonly used models e.g., of infinite coin tosses, some objections to
the use of Robinson's framework stem from hidden and unnecessary
assumptions that predetermine the choice of an Archimedean model.
Dropping such hidden assumptions enables alternative models of these
infinite processes, obtained via hyperfinite techniques.

Moreover, we have started addressing the claim by Pruss that
transferless extensions of the real numbers (such as the surreal
field, the Levi-Civita field, or the field of Laurent series) might be
more suitable for the development of infinitesimal probabilities.  The
proposal of working with such non-Archimedean fields ignores both the
Klein--Fraenkel criteria for gauging the applicability of theories
with infinitesimals, and the power and utility of the transfer
principle.  In \cite{BK20b} we show that, due to these limitations,
the measure theory on such transferless fields is less expressive than
the hyperfinite counting measures.

Pruss claims that ``[w]hatever you can do with hyperreals, you can do
with surreals'' and that ``[t]he fields of formal Laurent series and
the Levi--Civita field {\ldots} have the advantage of being elegantly
small''\;\cite{Pr18}.  However, in \cite{BK20b} we will see that, by
his own Theorem\;1, these fields suffer from underdetermination due to
the possibility of rescaling the infinitesimal part of the
probability, and do not possess a notion of internality that enables
one to escape such underdetermination in Robinson's framework.
Prussian Theorem\;$1$, when properly analyzed, boomerangs to undercut
his own underdetermination thesis.  Prussian Theorem\;$2$ is similarly
Boomerang\;$2$ due to the existence of nontrivial automorphisms for
all such transferless fields.%
\footnote{See Ehrlich (\cite{Eh94}, 1994, pp.\;253) for the existence
of nontrivial automorphisms of the surreals as an ordered field, and
Shamseddine (\cite{Sh11b}, 2011, p.\;224, Remark 3.17) for such
existence for all ordered extensions of~$\R$.}

Arguments based on the non-effectiveness of ultrafilters are not
limited to the work of Pruss; see e.g., Easwaran--Towsner \cite{et}.
In spite of an initial plausibility of such arguments against
Robinson's framework, the arguments dissolve upon closer inspection,
and even tend to prove the opposite of what their authors intended.
One can well wonder why such arguments \emph{ad\;ultrafiltrum} don't
succeed.  A recent development suggests a possible reason.  It turns
out that the main body of the applicable part of Robinson's framework
admits a formalisation that requires modest foundational means not
exceeding those required for traditional non-infinitesimal methods in
ordinary mathematics; see Hrbacek--Katz \cite{Hr20}.  Thus the alleged
non-effectiveness is simply not there to begin with.

\section*{Acknowledgments} 

We are grateful to Karel Hrbacek, Vladimir Kanovei, Karl Kuhlemann,
and David Sherry for insightful comments on earlier versions that
helped improve our article, and to anonymous referees for constructive
criticism.  The influence of Hilton Kramer (1928--2012) is obvious.

\end{document}